\documentclass{article}
\usepackage{amsfonts}
\usepackage{graphicx}
\usepackage{amsmath}
\usepackage{amssymb}

\newtheorem{theorem}{Theorem}
\newtheorem{acknowledgement}[theorem]{Acknowledgement}

\newtheorem{corollary}[theorem]{Corollary}

\newtheorem{example}[theorem]{Example}

\newenvironment{proof}[1][Proof]{\textbf{#1.} }{\ \rule{0.5em}{0.5em}}
\numberwithin{theorem}{section}
\input{tcilatex}
\usepackage{graphicx}

\begin{document}

\title{On Viviani's Theorem and its Extensions}
\author{Elias Abboud \\
{\small Beit Berl College, Doar Beit Berl, 44905 Israel }\\
{\small \ Email: eabboud @beitberl.ac.il}}
\maketitle

\begin{abstract}
Viviani's theorem states that the sum of distances from any point inside an
equilateral triangle to its sides is constant. We consider extensions of the
theorem and show that any convex polygon can be divided into parallel
segments such that the sum of the distances of the points to the sides on
each segment is constant. A polygon possesses the \textit{CVS }property if
the sum of the distances from any inner point to its sides is constant. An
amazing result, concerning the converse of Viviani's theorem is deduced;
Three non-collinear points which have equal sum of distances to the sides
inside a convex polygon, is sufficient for possessing the \textit{CVS }%
property. For concave polygons the situation is quite different, while for
polyhedra analogous results are deduced.

\qquad


\qquad

\end{abstract}

\section{\protect\bigskip Introduction}

Let $\mathcal{P}$ be a polygon or polyhedron, consisting of both boundary
and interior points. Define a distance sum function $\mathcal{V}:\mathcal{P}%
\rightarrow\mathbb{R},$ where for each point $P\in\mathcal{P}$ $\ $the value 
$\mathcal{V}(P)$ is defined as the sum of the distances from the point $P$
to the sides (faces) of $\mathcal{P}.$

We say that $\mathcal{P}$ has the \textit{constant Viviani sum }property,
abbreviated by the "\textit{CVS }property", if and only if the function $%
\mathcal{V}$ is constant.

Viviani (1622-1703), who was a student and assistant of Galileo, discovered
the theorem which states that equilateral triangles have the \textit{CVS }%
property. The theorem can be easily proved by an area argument; Joining a
point $P$ inside the triangle to its vertices divides it into three parts,
the sum of their areas will be equal to the area of the original one.
Therefore, $\mathcal{V}(P)$ will be equal to the height of the triangle and
the theorem follows. The importance of Viviani's theorem may be derived from
the fact that his teacher Torricelli (1608-1647) used it to locate the
Fermat point of a triangle \cite[pp. 443]{GT}.

Samelson \cite[pp. 225]{Sam} gave a proof of Viviaini's theorem that uses
vectors and Chen \& Liang \cite[pp. 390-391]{CL} used this vector method for
proving the converse of the theorem; If inside a triangle there is a
circular region in which $\mathcal{V}$ is constant then the triangle is
equilateral.

Kawasaki \cite[ pp. 213]{Kaw}, by a proof without words, uses only rotations
to establish Viviani's theorem. There is an extension of the theorem to all
regular polygons, by the area method: All regular polygons have the \textit{%
CVS }property. There is also an extension of the theorem to regular
polyhedrons, by a volume argument: All regular polyhedra have the \textit{%
CVS }property. Kawasaki, Yagi and Yanagawa \cite[pp. 283]{KYY} gave a
different proof for the regular tetrahedron.

What happens for general polygons and polyhedra?

Surprisingly, there is a strict correlation between Viviani's theorem, its
converse and extensions to linear programming.

This correlation is manifested by the following main result:

\begin{theorem}
\label{parallel segments}(a) Any convex polygon can be divided into parallel
segments such that $\mathcal{V}$ is constant on each segment. \ 

(b) Any convex polyhedron can be divided into parallel cross sections such
that $\mathcal{V}$ is constant on each cross section.
\end{theorem}

These segments or cross sections, on which $\mathcal{V}$ is constant, will
be called \textit{isosum layers (}or definitely\textit{\ isosum segments }and%
\textit{\ isosum cross sections)}. They are formed by the intersection of $%
\mathcal{P}$ and a suitable family of parallel lines (planes). The value of
the function $\mathcal{V}$ will increase when passing, in some direction,
from one isosum layer into another, unless $\mathcal{P}$ has the \textit{CVS}
property.

The correlation soon will be clear. Each linear programming problem is
composed of an \textit{objective function} and a \textit{feasible region }%
(see for example\textit{\ }\cite{LGR} or \cite{V}). Moreover, the objective
function divides the feasible region into \textit{isoprofit} layers, these
layers are parallel and consist of points on which the objective function
has constant value. Furthermore, moving in some direction will increase the
value of the objective function unless it is constant in the feasible region.

Because of this correlation we thus conclude the following amazing result,
concerning the converse of Viviani's theorem.

\begin{theorem}
\label{3 noncolinear}(a) If $\ \mathcal{V}$ takes equal values at three
non-collinear points, inside a convex polygon, then the polygon has the 
\textit{CVS }property.

(b) If $\mathcal{V}$ takes equal values at four non-coplanar points, inside
a convex polyhedron, then the polyhedron has the \textit{CVS }property.
\end{theorem}

The theorem tells us that measuring the distances from the sides of three
non-collinear points, inside a convex polygon, is sufficient for determining
if the polygon has the \textit{CVS }property. Likewise, measuring the
distances from the faces of four non-coplanar points, inside a convex
polyhedron, is sufficient for determining if the polyhedron possesses the 
\textit{CVS }property.

We then end with the following beautiful conclusions.

\begin{corollary}
\label{rotation symmetry polygon}(a) If there is an isometry of the plane
which fixes the polygon but not an isosum segment, then the polygon has the 
\textit{CVS }property.

(b) If a convex polygon possesses a \textit{rotational} symmetry, around a
central point, then the polygon has the \textit{CVS }property.

(c) If a convex polygon possesses a reflection symmetry across an axis $l,$
then the polygon has the CVS property or otherwise the isosum segments are
perpendicular to $l.$
\end{corollary}

While for polyhedra we have,

\begin{corollary}
\label{rotation symmetry polyhedron}(a) If there is an isometry of the space
which fixes the polyhedron but not an isosum cross section, then the
polyhedron has the \textit{CVS }property.

(b) If a convex polyhedron possesses two \textit{rotational} symmetries,
around different axes, then the polyhedron has the \textit{CVS }property.
\end{corollary}

From Corollary \ref{rotation symmetry polygon}, one can deduce that all
regular polygons have the \textit{CVS }property. Besides, any parallelogram
has this property, since it possesses a rotational symmetry around its
centroid by an angle of $180^{\circ}.$

Obviously, the existence of two reflection symmetries, by different axes, of
a polygon will imply a rotational symmetry and hence the polygon must own
the \textit{CVS }property.

Moreover, for triangles and quadrilaterals, the existence of a rotational
symmetry characterizes the possessing of \textit{CVS }property, since in
these cases the polygons would be only equilateral triangles and
parallelograms.

Consequently, an $n$-gon, for $n\geq5,$ that does not possess the \textit{%
CVS }property must have at most one symmetry which is the reflection
symmetry.

Analogously, by Corollary \ref{rotation symmetry polyhedron}, all regular
polyhedra and regular prisms have the \textit{CVS }property. Likewise, any
parallelepiped has the \textit{CVS }property. Since it possesses three
rotational symmetries by an angle of $180^{\circ},$ around the axes, each
passing through the centroids of a pair of parallel faces.

On the other hand, the property of possessing a symmetry does not
characterize all polygons (polyhedra) satisfying the \textit{CVS }property.
In section IV. we will validate the existence of a polygon with only one
reflection symmetry, an asymmetric polygon and a polyhedron with a
reflection symmetry only, which possess the \textit{CVS }property.

We will proceed as follows. In section II., we first introduce a linear
programming problem for general triangles. The main statement will be;

\textquotedblright A triangle has the \textit{CVS }property \textit{if and
only if} \ it is equilateral, \textit{if and only if} there are three
non-collinear points inside the triangle that have equal sum of distances
from the sides\textquotedblright.

Then we deal with general convex polygons and polyhedrons and the proof of
theorems (\ref{parallel segments}), (\ref{3 noncolinear}). Here we shall
rely on methods from analytic geometry because the use of coordinates. This
allows us to determine the line (plane) which the isosum layers are parallel
to.

In section III. we see what happens for concave polygons and polyhedra, then
in section IV. we compute some examples.

It is worthy mentioning that these results can be stated and generalized for 
$n-$dimensional geometry.

\section{Convex Polygons and Polyhedra}

\subsection{The case of triangle: linear programming approach}

Given a triangle $\triangle ABC$ \ let $a_{1},$ $a_{2},$ $a_{3}$ be the
lengths of the sides $BC,$ $AC$, $AB$ respectively. Let $P$ be a point
inside the triangle and let $h_{1},$ $h_{2},$ $h_{3}$ be the distances (the
lengths of the altitudes) of the point $P$ to the three sides respectively,
see Figure \ref{fig1}.

For $1\leq i\leq3,$ let $x_{i}=\frac{h_{i}}{\sum_{i=1}^{3}h_{i}},$ where as
previously defined, $\sum_{i=1}^{3}h_{i}=\mathcal{V}(P)$. Clearly, for each $%
1\leq i\leq3,$ we have $0\leq x_{i}\leq1$ and $\sum_{i=1}^{3}x_{i}=1.$
Denote $x=(x_{1},x_{2},x_{3})$ and consider the linear function in three
variables $F(x)=\sum_{i=1}^{3}a_{i}x_{i}.$ Now, this function is closely
related to the function $\mathcal{V}$ . Accurately, $F(x)=%
\sum_{i=1}^{3}a_{i}x_{i}=\frac {\sum_{i=1}^{3}a_{i}h_{i}}{\sum_{i=1}^{3}h_{i}%
}=\frac{2S}{\mathcal{V}(P)},$ where $S$ is the area of the triangle.

Consequently, $F(x)=\sum_{i=1}^{3}a_{i}x_{i}$ takes equal values in a subset
of points of the feasible region \textit{if and only if} the function $%
\mathcal{V}$ takes equal values at the corresponding points inside the
triangle.

Thus we may define the following linear programming problem;

The objective function is: 
\begin{equation*}
F(x)=\sum_{i=1}^{3}a_{i}x_{i}
\end{equation*}

subject to the following constraints:

\begin{equation*}
\left\{ 
\begin{array}{c}
\sum_{i=1}^{3}x_{i}\leq1\qquad \\ 
x_{i}\geq0,\text{\qquad}1\leq i\leq3%
\end{array}
.\right.
\end{equation*}

Now, solving the problem means maximizing or minimizing the objective
function in the feasible region, and this optimal value must occur at some
corner point. For that, one may use the simplex method, which is simple in
this case because the simplex tableau contains only two rows. But rather
than using this algebraic method, and since we are not seeking for optimal
values, and for better understanding the correlation with Viviani's theorem,
we use a geometric method. The feasible region will be the right pyramid
with vertex $(0,0,1)$ and basis vertices $(0,0,0),(1,0,0),(0,1,0)$ and the
isoprofit planes\ are obtained by taking $F(x)=\sum_{i=1}^{3}a_{i}x_{i}=c,$
where $c$ is constant, see Figure \ref{fig2}.

These isoprofit planes\ will meet the face $\sum_{i=1}^{3}x_{i}=1$ across a
line which in its turn will correspond to an isosum segment inside the
triangle $\triangle ABC.$ The value of the function $\mathcal{V}$ will be
constant for all points on the isosum segment. Moreover, moving in a
perpendicular direction to the isosum segments will increase the value of $%
\mathcal{V}$. Hence, $\mathcal{V}$ does not take the same value at any three
non-collinear points inside the triangle.

This is true in general, except for one case where the isoprofit planes,\ $%
\sum_{i=1}^{3}a_{i}x_{i}=c,$ are parallel to the face $\sum
_{i=1}^{3}x_{i}=1. $ This exceptional case occurs \textit{if and only if}
the vectors $(a_{1},a_{2},a_{3})$ and $(1,1,1)$ are linearly dependent i.e., 
\begin{equation*}
(a_{1},a_{2},a_{3})=\lambda(1,1,1),
\end{equation*}

and this happens \textit{if and only if} the triangle is equilateral. Thus
proving the following theorem.

\begin{theorem}
(a) Any triangle can be divided into parallel segments such that $\mathcal{V}
$ is constant on each segment.

(b) The following conditions are equivalent

\begin{itemize}
\item The\textit{\ triangle }$\triangle ABC$ \textit{has the CVS property. }

\item \textit{There are three non-collinear points, inside the triangle, at
which }$\mathcal{V}$ takes the same value\textit{.}

\item \textit{\ }$\triangle ABC$ is \textit{equilateral. }
\end{itemize}
\end{theorem}

\subsection{Extension for convex polygons: analytic geometry approach}

Inspired by the geometric method of the linear programming problem defined
in the previous subsection, we attack the general case by using analytic
geometry techniques.

Firstly, we aim to prove Theorems \ref{parallel segments}(a) and \ref{3
noncolinear}(b).

Given a polygon with $n$ sides, we embed it in a Cartesian plane. Suppose
that its sides lie on lines with equations; 
\begin{equation}
\alpha_{i}X+\beta_{i}Y+\gamma_{i}=0.  \label{0}
\end{equation}

Since the convex polygon lies on the same side of line (\ref{0}) then the
expression $\alpha_{i}x+\beta_{i}y+\gamma_{i}$ has unchanging sign, for all $%
P=(x,y)$ inside the polygon.

Thus the distance, $h_{i},$ of the point $P$ from each side of the polygon
is given by the equation; 
\begin{equation*}
h_{i}=(-1)^{\varepsilon_{i}}\frac{\alpha_{i}x+\beta_{i}y+\gamma_{i}}{\sqrt{%
\alpha_{i}^{2}+\beta_{i}^{2}}},
\end{equation*}

where $\varepsilon_{i}\in\{0,1\}.$

Therefore, the function $\mathcal{V}$ is given by a linear expression 
\begin{equation}
\mathcal{V}(x,y)\mathcal{=}\sum_{i=1}^{n}(-1)^{\varepsilon_{i}}\frac {%
\alpha_{i}x+\beta_{i}y+\gamma_{i}}{\sqrt{\alpha_{i}^{2}+\beta_{i}^{2}}}.
\label{1}
\end{equation}

Letting this sum equals some constant $c,$ we get the equation of the line
which the isosum segments are parallel to; 
\begin{equation}
\overset{n}{\underset{i=1}{\sum}}(-1)^{\varepsilon_{i}}\frac{%
\alpha_{i}x+\beta_{i}y+\gamma_{i}}{\sqrt{\alpha_{i}^{2}+\beta_{i}^{2}}}=c.
\label{2}
\end{equation}

If the point $P=(x,y)$ is restricted to be inside the polygon then for
different values of $c,$ we get parallel segments. On each such segment, the
function $\mathcal{V}$ takes the constant value $c$. This proves Theorem \ref%
{parallel segments} (a).

The equation given by (\ref{1}) is independent of $(x,y)$ \textit{if and
only if} the variable part of \ the function $\mathcal{V}$ vanishes. But
then $\mathcal{V}$ is constant and so the polygon has the \textit{CVS }%
property.

If there are three non-collinear points inside the polygon at which $%
\mathcal{V}$ takes the same value then there exist two different isosum
segments on which the function $\mathcal{V}$ takes the same value. This
happens \textit{if and only if} \ $\mathcal{V}$ is constant and so the
polygon once again has the \textit{CVS }property. Hence Theorem \ref{3
noncolinear}(a) is proved.

We turn now to proving Corollary \ref{rotation symmetry polygon}.

(a) If there is an isometry of the plane which fixes the polygon but not
some isosum layer, then this will assure the existence of three
non-collinear points inside the polygon at which $\mathcal{V}$ takes equal
values. That's because an isometry preserves distances and the sets of
boundary and inner points. Thus Theorem \ref{3 noncolinear}(a) yields that
the polygon has the \textit{CVS }property.

(b) Follows from part (a).

(c) If the isosum segments are not fixed by the reflection, then according
to part (a) the polygon has the \textit{CVS} property. Otherwise, if fixed,
then the isosum segments must be perpendicular to the reflection axis.

\subsection{Extension for convex polyhedra}

The proof for a Polyhedron is similar up to minor modifications. The faces
lie on planes with equations; 
\begin{equation*}
\alpha_{i}x+\beta_{i}y+\gamma_{i}z+\delta_{i}=0,
\end{equation*}

and the linear function $\mathcal{V}$ becomes; 
\begin{equation*}
\mathcal{V}(x,y,z)=\overset{n}{\underset{i=1}{\sum}(-1)^{\varepsilon_{i}}}%
\frac{\alpha_{i}x+\beta_{i}y+\gamma_{i}z+\delta_{i}}{\sqrt{%
\alpha_{i}^{2}+\beta_{i}^{2}+\gamma_{i}^{2}}},
\end{equation*}

where $\varepsilon_{i}\in\{0,1\}.$

The same argument yields the result for the polyhedron (part (b) of Theorems %
\ref{parallel segments} and \ref{3 noncolinear}).

The proof of Corollary \ref{rotation symmetry polyhedron} is similar to
corollary \ref{rotation symmetry polygon};

(a) Follows from part (b) of Theorem \ref{3 noncolinear}.

(b) Two rotational symmetries of a polyhedron around different axes,
guarantee the existence of four non-coplanar points inside the polyhedron at
which $\mathcal{V}$ takes equal values. By Theorem \ref{3 noncolinear} (b),
the result follows.

\section{Concave Polygons and Polyhedra}

The situation for concave polygons and polyhedra is quite different.
Theorems \ref{parallel segments} and \ref{3 noncolinear} are no longer
valid. Moreover, concave polygons and polyhedra don't have the \textit{CVS}
property. Surprisingly, with a little more elaborate effort, one might
rather find a generality of Theorem \ref{parallel segments}. Concerning
polygons, the crucial point which makes the difference is that the points
inside a concave polygon are no longer lie on the same side of each boundary
line. This was a key for defining the distance sum function $\mathcal{V}$.

We turn to giving an example to illustrate the theme. Let $ABCD$ be the
concave kite with vertices $(0,8),(-6,0),(0,2.5),(6,0)$ respectively. Let $%
l_{1},l_{2},l_{3},l_{4}$ be the lines containing the sides $AB,BC,CD,DA$ and 
$E,F$ be the intersection points of $l_{1},l_{3}$ and $l_{2},l_{4}$
respectively, see Figure \ref{fig3}. Then the concave kite $ABCD$ is divided
into three distinct convex polygonal regions, namely; $AECF,EBC$ and $FCD.$
All points inside any region lie on the same side of each $l_{i},$ $1\leq
i\leq4.$ To explain that, note that any line $l$ divides the plane into two
half-planes $"O_{l}"$ which contains the origin and its complement $%
"O_{l}^{c}"$ . The following table shows the location of the points inside
any region relative to the lines $l_{i},$ $1\leq i\leq4.$

\begin{equation*}
\begin{array}{cc}
& \text{location relative to} \\ 
\begin{array}{c}
\text{region} \\ 
AECF \\ 
EBC \\ 
FCD%
\end{array}
& 
\begin{array}{cccc}
l_{1} & l_{2} & l_{3} & l_{4} \\ 
O_{l_{1}} & O_{l_{2}}^{c} & O_{l_{3}}^{c} & O_{l_{4}} \\ 
O_{l_{1}} & O_{l_{2}}^{c} & O_{l_{3}} & O_{l_{4}} \\ 
O_{l_{1}} & O_{l_{2}} & O_{l_{3}}^{c} & O_{l_{4}}%
\end{array}%
\end{array}%
\end{equation*}

Note that any two neighboring regions, having a common edge, differ only in
one entry in this table. For instance, the points of $AECF$ and $EBC$ lie on
the same side of $l_{1},l_{2},l_{4}$ and on opposite sides of $l_{3}.$ The
implementation is that the distance sum function $\mathcal{V}$ will be a
split function composed of three components;%
\begin{equation*}
\mathcal{V}_{ABCD}(P)=\left\{ 
\begin{array}{c}
\mathcal{V}_{AECF}(P),\text{ \ \ \ \ }P\text{ inside }AECF \\ 
\mathcal{V}_{EBC}(P),\text{ \ \ \ \ \ \ \ }P\text{ inside }EBC \\ 
\mathcal{V}_{FCD}(P),\text{ \ \ \ \ \ \ \ }P\text{ inside }FCD%
\end{array}
,\right.
\end{equation*}

where each component is given by a linear expression as in equation (\ref{1}%
), definitely;

\begin{align*}
\mathcal{V}_{AECF}(P) & =-\frac{-8x+6y-48}{10}+\frac{-5x+12y-30}{13}+\frac{%
5x+12y-30}{13} \\
& -\frac{8x+6y-48}{10} \\
\mathcal{V}_{EBC}(P) & =-\frac{-8x+6y-48}{10}+\frac{-5x+12y-30}{13}-\frac{%
5x+12y-30}{13} \\
& -\frac{8x+6y-48}{10} \\
\mathcal{V}_{DCF}(P) & =-\frac{-8x+6y-48}{10}-\frac{-5x+12y-30}{13}+\frac{%
5x+12y-30}{13} \\
& -\frac{8x+6y-48}{10}
\end{align*}

This says that each region can be divided into parallel isosum segments but
in a different direction. In our example, the isosum segments of $AECF$ are
parallel to a line of the form $y=c,$ the isosum segments of $EBC$ are
parallel to the line $100x+156y=0$ and the isosum segments of $DCF$ are
parallel to the line $100x-156y=0.$ Moreover, two isosum segments from
neighboring regions which meet at a point on the common edge, define three
non-collinear points with the same distance sum from the sides of $ABCD,$
see Figure \ref{fig3}.

These ideas can be generalized to any concave polygon or polyhedron. Thus we
have the following theorem;

\begin{theorem}
(a) Any concave polygon can be divided into convex polygonal regions such
that each region can be divided into parallel isosum segments. Moreover,
isosum segments of neighboring regions have different directions.

(b) There are three non-collinear points inside a concave polygon which have
equal distance sum from the sides.

(c) Any concave polyhedron can be divided into convex polyhedral regions
such that each region can be divided into parallel isosum cross sections.
Moreover, isosum cross sections of neighboring regions have different
directions.

(d) There are four non-coplanar points inside a concave polyhedron which
have equal distance sum from the faces.

(e) Concave polygons and polyhedra do not possess the \textit{CVS} property.

\begin{proof}
We shall outline the proof for polygons.

(a) Let $\mathcal{P}$ a concave polygon. Extend the sides and construct all
possible intersection points in $\mathcal{P}$ of the boundary lines, see
Figure \ref{fig4}. In this way one gets a partition of $\mathcal{P}$\ into $%
m $ convex polygonal regions $\mathcal{P}_{1},\mathcal{P}_{2},...,\mathcal{P}%
_{m}.$ Note that each side of any $\mathcal{P}_{j}$ lie on a boundary line
of $\mathcal{P}$. The points inside any two neighboring regions $\mathcal{P}%
_{i},\mathcal{P}_{j}$ with a common edge $LM$, lie on the same side of each
boundary line of $\mathcal{P}$ except one, namely, the boundary line $l_{LM}$
containing $LM,$ see Figure \ref{fig4}$.$ We claim that the conjugation of
these two neighboring regions $\mathcal{P}_{i}\cup\mathcal{P}_{j}$, will
remain a convex polygonal region. This is true since the two sides of $%
\mathcal{P}_{i}$ and $\mathcal{P}_{j}$ which meet at $L$ (or $M)$ lie on the
same boundary line of $\mathcal{P}$, by the construction (in Figure \ref%
{fig4} the points $K,L,D$ are collinear, likewise $I,M,D$ are collinear).
Thus all points of $\mathcal{P}_{i}\cup\mathcal{P}_{j}$ lie on the same side
of the boundary lines of $\mathcal{P}$ except $l_{LM}.$

Now, the distance sum function $\mathcal{V}$ \ of $\mathcal{P}$ is defined
as a split function;%
\begin{equation*}
\mathcal{V}_{\mathcal{P}}(P)=\left\{ 
\begin{array}{c}
\mathcal{V}_{\mathcal{P}_{1}}(P),\text{ \ \ \ \ }P\text{ inside }\mathcal{P}%
_{1} \\ 
\mathcal{V}_{\mathcal{P}_{2}}(P),\text{ \ \ \ \ }P\text{ inside }\mathcal{P}%
_{2} \\ 
. \\ 
. \\ 
\mathcal{V}_{\mathcal{P}_{m}}(P),\text{ \ \ \ \ }P\text{ inside }\mathcal{P}%
_{m}%
\end{array}
\right. .
\end{equation*}
The points on a common edge of neighboring regions can be attached to any
one of both. Each convex polygonal region has parallel isosum segments
according to the linear expression which defines $\mathcal{V}_{\mathcal{P}%
_{i}},1\leq i\leq m.$ Isosum segments of neighboring regions $\mathcal{P}%
_{i},\mathcal{P}_{j}$ have different directions since the linear expressions
of $\mathcal{V}_{\mathcal{P}_{i}},\mathcal{V}_{\mathcal{P}j}$ differ only in
one sign as explained above.

(b) Take any point $P$ on a common edge of two neighboring regions $\mathcal{%
P}_{i},\mathcal{P}_{j}$ and let $s_{1},s_{2}$ be two isosum segments of $%
\mathcal{P}_{i},\mathcal{P}_{j}$ respectively, issuing from $P.$ Let $%
Q_{1},Q_{2}$ two points on $s_{1},s_{2}$ respectively. Since $s_{1},s_{2}$
have different directions then $P,$ $Q_{1},Q_{2}$ are non-collinear and have
the same distance sum from the sides.

(e) The distance sum functions $\mathcal{V}_{\mathcal{P}_{i}},\mathcal{V}_{%
\mathcal{P}j}$ of two neighboring regions $\mathcal{P}_{i},\mathcal{P}_{j}$
have the same linear expressions with the same signs except one; $\mathcal{V}%
_{\mathcal{P}_{i}}=\underset{i=1}{\overset{m}{\sum}}%
(a_{k}x+b_{k}y+c_{k})+(a_{0}x+b_{0}y+c_{0})$ and $\mathcal{V}_{\mathcal{P}j}=%
\underset{i=1}{\overset{m}{\sum}}(a_{k}x+b_{k}y+c_{k})-(a_{0}x+b_{0}y+c_{0}) 
$. Thus $\mathcal{V}_{\mathcal{P}}$ cannot be constant on $\mathcal{P}$.
\end{proof}
\end{theorem}

\section{Examples}

In the following we compute the equations of the isosum layers for
particular polygons and polyhedra. It is an easy matter to check the
results, by constructing the figures using a computer program such as the
geometer sketchpad or wingeom.

\begin{example}
(Kite)
\end{example}

Construct a kite in the Cartesian plane with vertices, $A=(0,\beta
),B=(-\alpha,0),C=(0,\gamma),D=(\alpha,0),$ where, $\alpha,\beta>0$ and $%
\gamma<0.$

The equations of the lines containing the sides $AB,BC,CD,DA,$ respectively
are;

\begin{align*}
-\alpha y+\beta x+\alpha\beta & =0 \\
-\alpha y+\gamma x+\alpha\gamma & =0 \\
\alpha y+\gamma x-\alpha\gamma & =0 \\
\alpha y+\beta x-\alpha\beta & =0.
\end{align*}

Thus the function $\mathcal{V}$ is given by;

\begin{equation*}
\mathcal{V}_{ABCD}=\frac{-\alpha y+\beta x+\alpha\beta}{\sqrt{%
\alpha^{2}+\beta^{2}}}-\frac{-\alpha y+\gamma x+\alpha\gamma}{\sqrt{%
\alpha^{2}+\gamma^{2}}}+\frac{\alpha y+\gamma x-\alpha\gamma}{\sqrt{%
\alpha^{2}+\gamma^{2}}}-\frac{\alpha y+\beta x-\alpha\beta}{\sqrt{%
\alpha^{2}+\beta^{2}}}.
\end{equation*}

Equivalently, 
\begin{equation}
\mathcal{V}_{ABCD}=\frac{-2\alpha y}{\sqrt{\alpha^{2}+\beta^{2}}}+\frac{%
2\alpha\beta}{\sqrt{\alpha^{2}+\beta^{2}}}+\frac{2\alpha y}{\sqrt{%
\alpha^{2}+\gamma^{2}}}-\frac{2\alpha\gamma}{\sqrt{\alpha^{2}+\gamma^{2}}}.
\label{3}
\end{equation}

As a result we see that the isosum segments have simple equations, namely, $%
y=c.$

Obviously,\ equation (\ref{3}), shows that the $\mathcal{V}_{ABCD}$ is
independent of variables \textit{if and only if} \ $\beta=\gamma.$ This
means that the polygon is in fact a parallelogram and hence, it owns the 
\textit{CVS}\ property. On the other side, if $\beta\neq\gamma,$ then the
kite does not own the \textit{CVS}\ property, but it can be divided into
isosum segments which are parallel to the diagonal $BD$, which is consistent
with Corollary \ref{rotation symmetry polygon} (c).

\begin{example}
(Isosceles\ triangle)
\end{example}

In the previous example if one considers $P$ as a point inside the isosceles
triangle $ABD$ then, 
\begin{equation*}
\mathcal{V}_{ABD}=\frac{-\alpha y+\beta x+\alpha\beta}{\sqrt{%
\alpha^{2}+\beta^{2}}}-\frac{\alpha y+\beta x-\alpha\beta}{\sqrt{%
\alpha^{2}+\beta^{2}}}+y.
\end{equation*}

Equivalently, 
\begin{equation*}
\mathcal{V}_{ABD}=\frac{(\sqrt{\alpha^{2}+\beta^{2}}-2\alpha)y}{\sqrt {%
\alpha^{2}+\beta^{2}}}+\frac{2\alpha\beta}{\sqrt{\alpha^{2}+\beta^{2}}}.
\end{equation*}

The function $\mathcal{V}_{ABD}$ is independent of variables \textit{if and
only if} \ $\beta=\sqrt{3}\alpha.$ In which case the triangle is in fact
equilateral and hence it has the \textit{CVS}\ property.

If $\beta\neq\sqrt{3}\alpha,$ then it can be divided into isosum segments
parallel to the base, which is consistent with Corollary \ref{rotation
symmetry polygon} (c)

\begin{example}
(Quadrilateral)
\end{example}

Let be given the quadrilateral with vertices $(0,0),(3,0),(1,2),(0,1).$ \
Direct computations will show that the isosum segments are parallel to the
line $y=(1+\sqrt{2})x,$ \ see Figure \ref{fig5}. This quadrilateral does not
have the \textit{CVS }property.

\begin{example}
\label{pentagon}(Pentagon with a reflection symmetry)
\end{example}

The pentagon $ABCDE$ with vertices $A=(0,3+\sqrt{3}%
),B=(-1,3),C=(-1,0),D=(1,0),E=(1,3)$ and sides $2,3,2,3,2$ respectively,
owns the \textit{CVS }property although it has one reflection symmetry and
no rotational symmetries. It is a straightforward matter, to show the
possessing of the \textit{CVS }property, since it is composed of a rectangle
and an equilateral triangle.

In general, suppose that $ABCDE$ is a pentagon where $ABE$ is an equilateral
triangle with side-length $a$ and height $h,$ and $BCDE$ is a rectangle with
side-lengths $a,b.$ Then the distance sum of an inner point from its sides
is $a+b+h.$ This is obvious if the point is inside the triangle. If the
point is inside the rectangle then observe that it is located on the base of
another equilateral triangle obtained by extending the sides $AB,AE$ and the
result easily follows.

\begin{example}
\label{nosymmetries}(Asymmetric pentagon)
\end{example}

Construct any pentagon $ABCDE$ with angles $70^{\circ},110^{\circ},130^{%
\circ },60^{\circ}$ and $170^{\circ}$. Here again it is a straightforward
matter, to show that $ABCDE$ possess the \textit{CVS }property, since it is
composed of a parallelogram and an equilateral triangle. This pentagon has
no symmetries.

Observe that in this case, the distance sum of an inner point from the sides
is $a+b+h,$ where $a,b$ are the distances between the opposite sides of the
parallelogram and $h$ is the height of the equilateral triangle.

\begin{example}
(Equiangular polygon)
\end{example}

Any equiangular polygon $\mathcal{P}$ with $n$ sides has the \textit{CVS}
property. This can be demonstrated by the following argument; Locate inside
the equiangular polygon a regular $n$-gon, $\mathcal{P}_{r}$. Rotate $%
\mathcal{P}_{r}$ around its centroid until one of its sides gets parallel to
one side of $\mathcal{P}$. But then all corresponding sides of both polygons
get parallel. Let $\mathcal{V}_{\mathcal{P}},$ $\mathcal{V}_{\mathcal{P}%
_{r}} $ be the distance sum functions defined on $\mathcal{P}$ and $\mathcal{%
P}_{r} $ respectively. Then for any point $P$ inside $\mathcal{P}_{r}$ we
have,

\begin{equation*}
\mathcal{V}_{\mathcal{P}}(P)=\mathcal{V}_{\mathcal{P}_{r}}(P)+c,
\end{equation*}

where $c$ represents the sum of distances between the parallel sides of $%
\mathcal{P}$ and $\mathcal{P}_{r}.$

Since $\mathcal{V}_{\mathcal{P}_{r}}$ is constant inside $\mathcal{P}_{r}$
and also $c$ is constant then $\mathcal{V}_{\mathcal{P}}$ is constant inside 
$\mathcal{P}_{r}.$ By Theorem \ref{3 noncolinear}(a), $\mathcal{V}_{\mathcal{%
P}}$ is constant on $\mathcal{P}$.

\begin{example}
\label{polyhedronrefsym}(Polyhedron with a reflection symmetry)
\end{example}

Construct the right prism for which the pentagon $ABCDE$ from the previous
example congruent to its bases. This prism has only one symmetry, namely,
the reflection symmetry by the plane parallel to its two bases and passes in
the middle. This polyhedron has the \textit{CVS }property. Since vertically
the sum of the two altitudes to the bases from any inner point is constant
and equals to the height of the prism and, horizontally, the sum of the
altitudes to the faces is also constant by the previous example.

\begin{example}
\label{pyramidrotsym}(Pyramid with a rotational symmetry around one axis)
\end{example}

Construct in the space, the rectangular pyramid with vertices $A=(0,0,\alpha
)$, $B=(\beta,0,0),C=(0,\gamma,0),D=(-\beta,0,0)$ and $E=(0,-\gamma,0),$
where $\alpha,\beta,\gamma>0.$ The faces $ABC,ACD,ADE$ and $AEB$ lie on the
following planes respectively; 
\begin{align*}
\alpha\gamma x+\alpha\beta y+\beta\gamma z & =\alpha\beta\gamma \\
-\alpha\gamma x+\alpha\beta y+\beta\gamma z & =\alpha\beta\gamma \\
-\alpha\gamma x-\alpha\beta y+\beta\gamma z & =\alpha\beta\gamma \\
\alpha\gamma x-\alpha\beta y+\beta\gamma z & =\alpha\beta\gamma.
\end{align*}

Letting $\Delta=\sqrt{\alpha^{2}\gamma^{2}+\alpha^{2}\beta^{2}+\beta^{2}%
\gamma^{2}}$ then the altitudes from the point $P=(x,y,z)$ to the faces are
given by; 
\begin{align*}
h_{BCDE} & =z \\
h_{ABC} & =\frac{-1}{\Delta}(\alpha\gamma x+\alpha\beta y+\beta\gamma
z-\alpha\beta\gamma) \\
h_{ACD} & =\frac{-1}{\Delta}(-\alpha\gamma x+\alpha\beta y+\beta\gamma
z-\alpha\beta\gamma) \\
h_{ADE} & =\frac{-1}{\Delta}(-\alpha\gamma x-\alpha\beta y+\beta\gamma
z-\alpha\beta\gamma) \\
h_{AEB} & =\frac{-1}{\Delta}(\alpha\gamma x-\alpha\beta y+\beta\gamma
z-\alpha\beta\gamma)
\end{align*}

Thus 
\begin{equation*}
\mathcal{V}_{ABCDE}=(1-\frac{4\beta\gamma}{\Delta})z+\frac{%
4\alpha\beta\gamma }{\Delta}.
\end{equation*}

Evidently, the isosum layers are parallel to the base of the pyramid and it
has the \textit{CVS} property if and only if 
\begin{equation*}
\alpha=\frac{\sqrt{15}\beta\gamma}{\sqrt{\beta^{2}+\gamma^{2}}}.
\end{equation*}

Figure \ref{fig6} shows a special case where, $\beta=\gamma=1$ and $\alpha=%
\sqrt{\frac{15}{2}}.$

\section{Concluding Remarks}

\begin{itemize}
\item Corollary \ref{rotation symmetry polygon} part (c) can be extended for
polyhedra but the property of being perpendicular to the reflection plane
would not uniquely determine the isosum cross sections.

\item It is possible to state an algebraic necessary and sufficient
condition for the \textit{CVS }property, using the expression of $\mathcal{V}
$ given in (\ref{1}). But a geometric one is more favorable.

Thus the substantial question is how can one characterizes, geometrically,
all those polygons and polyhedra that satisfy the\textit{\ CVS} property?
\end{itemize}

\begin{acknowledgement}
The author is indebted for all colleagues who took part in reviewing the
manuscript. Special thanks are due to the anonymous referees for their
valuable comments that improved the exposition of the paper.
\end{acknowledgement}

\bigskip

\FRAME{ftbphFU}{2.7086in}{2.1335in}{0pt}{\Qcb{ A point $P$ inside the
triangle with \ distances $h_{1},$ $h_{2},$ $h_{3}$.}}{\Qlb{fig1}}{%
triangle.eps}{\special{language "Scientific Word";type
"GRAPHIC";maintain-aspect-ratio TRUE;display "USEDEF";valid_file "F";width
2.7086in;height 2.1335in;depth 0pt;original-width 7.7012in;original-height
11.1431in;cropleft "0.2606";croptop "0.6394";cropright "0.6087";cropbottom
"0.4505";filename
'../../Elias/epsfiles/VivianiTheorem/triangle.eps';file-properties "XNPEU";}}

\bigskip

\FRAME{ftbpFU}{4.7167in}{3.1393in}{0pt}{\Qcb{The feasible region with an
isoprofit plane: $F(x)=0.$}}{\Qlb{fig2}}{feasible.eps}{\special{language
"Scientific Word";type "GRAPHIC";maintain-aspect-ratio TRUE;display
"USEDEF";valid_file "F";width 4.7167in;height 3.1393in;depth
0pt;original-width 7.7012in;original-height 11.1431in;cropleft
"0.1301";croptop "0.6846";cropright "0.7390";cropbottom "0.4054";filename
'../../Elias/epsfiles/feasible.eps';file-properties "XNPEU";}}

\FRAME{ftbpFU}{4.7167in}{3.1419in}{0pt}{\Qcb{The isosum segments have
different directions and $X,Y,Z$ are non-collinear points with equal
distance sum from the sides of $ABCD.$}}{\Qlb{fig3}}{concavepolygon.eps}{%
\special{language "Scientific Word";type "GRAPHIC";maintain-aspect-ratio
TRUE;display "USEDEF";valid_file "F";width 4.7167in;height 3.1419in;depth
0pt;original-width 7.7012in;original-height 11.1431in;cropleft
"0.1955";croptop "0.6397";cropright "0.8044";cropbottom "0.3602";filename
'../../Elias/epsfiles/concavepolygon.eps';file-properties "XNPEU";}}

\FRAME{ftbpFU}{3.2128in}{3.1401in}{0pt}{\Qcb{Partitioning the concave
polygon $ABCDEFGH$ into convex polygonal regions.}}{\Qlb{fig4}}{conp.eps}{%
\special{language "Scientific Word";type "GRAPHIC";maintain-aspect-ratio
TRUE;display "USEDEF";valid_file "F";width 3.2128in;height 3.1401in;depth
0pt;original-width 7.7012in;original-height 11.1431in;cropleft
"0.2603";croptop "0.8196";cropright "0.6740";cropbottom "0.5403";filename
'../../Elias/epsfiles/conp.eps';file-properties "XNPEU";}}

\FRAME{ftbpFU}{2.7138in}{2.1344in}{0pt}{\Qcb{The isosum segments are
parallel to the line $y=(1+\protect\sqrt{2})x.$}}{\Qlb{fig5}}{%
quadrilateral.eps}{\special{language "Scientific Word";type
"GRAPHIC";maintain-aspect-ratio TRUE;display "USEDEF";valid_file "F";width
2.7138in;height 2.1344in;depth 0pt;original-width 7.7012in;original-height
11.1431in;cropleft "0.3907";croptop "0.6393";cropright "0.7394";cropbottom
"0.4503";filename
'../../Elias/epsfiles/VivianiTheorem/quadrilateral.eps';file-properties
"XNPEU";}}

\FRAME{ftbpFU}{2.7051in}{3.1401in}{0pt}{\Qcb{A square pyramid which has the 
\textit{CVS} property. \ }}{\Qlb{fig6}}{pyramid.eps}{\special{language
"Scientific Word";type "GRAPHIC";maintain-aspect-ratio TRUE;display
"USEDEF";valid_file "F";width 2.7051in;height 3.1401in;depth
0pt;original-width 7.7012in;original-height 11.1431in;cropleft
"0.3259";croptop "0.6394";cropright "0.6738";cropbottom "0.3600";filename
'../../Elias/epsfiles/VivianiTheorem/pyramid.eps';file-properties "XNPEU";}}

\end{document}